\documentclass[12pt]{article}
\usepackage{amsmath}
\usepackage{amssymb}
\usepackage{amsfonts}

\newtheorem{thm}{Theorem}
\newtheorem{prop}[thm]{Proposition}
\newtheorem{lemma}[thm]{Lemma}

\def\pa{\partial}

\title{Complete coalescent diagram of the Painlev\'e equations}
\author{Yousuke Ohyama \\
\small  Graduate School of Information Science and Technology,\\
\small Osaka University\\ {}\\
Shoji Okumura\\
\small  Graduate School of Science, Osaka University\\
}

\begin{document}
\maketitle
\begin{abstract} We will revise Garnier-Okamoto's coalescent diagram of
isomonodromic deformations and give a complete coalescent diagram. In our 
viewpoint, we have ten types of isomonodromic deformations and two of them 
give the same type of the Painlev\'e equation. We can naturally put the thirty-fourth 
Painlev\'e equation in our  diagram, which corresponds to the Flaschka-Newell form of the second 
Painlev\'e equation.
\end{abstract}

\section{Introduction}
In this paper, we will revise Garnier-Okamoto's coalescent diagram of
isomonodromic deformations \cite{Oka} and will show a  complete coalescent diagram. 
In the original form, the  Painlev\'e equations are classified in six types. But in our 
picture, there exist ten different types of isomonodromic deformations.
Since two of them give the same type of the Painlev\'e equation, the Painlev\'e equations 
are classified in eight different types. 
 
We will also show that  the Painlev\'e 
equations are classified in five types as   a nonlinear single equation.
Especially we show a unified equation of the fourth Painlev\'e equation and 
the thirty-fourth Painlev\'e equation. These five types are classified into fourteen types 
by scaling transformations. 
We exclude four types of them since they are quadrature. The remaining  ten types of equations
correspond to the different singularity types of isomonodromic deformations. In our form, 
it is easy to understand the relation between the type of the Painlev\'e equations 
and the singularity type of isomonodromic deformations.

It is known that different forms of isomonodromic deformations 
\begin{eqnarray*}
     \frac{\pa Y}{\pa x} =A (x,t)Y,\quad  \frac{\pa Y}{\pa t} =B (x,t)Y  \label{deform}
\end{eqnarray*} 
exist for some types of the Painlev\'e equations. 
 One of the most famous example is
the Flaschka-Newell form \cite{FN} and the Miwa-Jimbo form \cite{JM2} for the second 
Painlev\'e equation P2($\alpha$)
\begin{equation}\label{0:p2}
y^{\prime\prime} = 2y^3 + t y +\alpha.
\end{equation}
The Flaschka-Newell form (FN) is 
\begin{equation}\label{FN}\begin{aligned}
A^{FN}(x,t)&=
-4 \begin{pmatrix}  x^2 &  y x\\   y x &   -x^2\end{pmatrix}+
\begin{pmatrix}  t+2y^2 &  -2z \\  2z &  -t-2y^2 \end{pmatrix}
- \begin{pmatrix}0 &  \alpha \\   \alpha & 0\end{pmatrix}\frac1x,\\
B^{FN}(x,t)&=  \begin{pmatrix} 1 & 0 \\ 0 & -1 \end{pmatrix}x +
 \begin{pmatrix} 0 & y \\  y & 0 \end{pmatrix}.
\end{aligned}\end{equation}
The Miwa-Jimbo form (MJ) is 
\begin{equation}\label{MJ}\begin{aligned}
     A^{MJ}(x,t)&=  \begin{pmatrix} 1 &0 \\ 0 &-1
                   \end{pmatrix}
             x^2+
            \begin{pmatrix} 0 &u \\ \frac{-2}{u}z &0
                   \end{pmatrix}x+
            \begin{pmatrix} z+\frac{t}{2} &-uy \\
                            \frac{-2}{u}(\theta+yz) &-z-\frac{t}{2}
                   \end{pmatrix},\\
     B^{MJ}(x,t)&=\frac{x}{2} \begin{pmatrix} 1 &0 \\ 0 &-1
                               \end{pmatrix}+
                        \frac{1}{2} 
                        \begin{pmatrix} 0 &u \\ -\frac{2}{u}z &0
                               \end{pmatrix}.
\end{aligned}\end{equation}
%% \begin{equation}\begin{aligned}
%% A^{FN}(x,t)&=
%% 4 \begin{pmatrix}-i x^2 & y x\\  y x & i x^2\end{pmatrix}+
%% \begin{pmatrix} -it-2iy^2 &  2iz \\  -2iz &  it+2iy^2 \end{pmatrix}
%% - \begin{pmatrix}0 &  \alpha \\   \alpha & 0\end{pmatrix}\frac1x,\\
%% B^{FN}(x,t)&=-i \begin{pmatrix} 1 & 0 \\ 0 & -1 \end{pmatrix}x +
%%  \begin{pmatrix} 0 & y \\  y & 0 \end{pmatrix}.
%% \end{aligned}\end{equation}
We take a slightly different form from the original
Flaschka-Newell form. Our form is a `real' form  in a sense.
$A^{FN}(x,t)$ has an irregular singularity 
of the Poincar\'e rank three at $x=\infty$ and  a regular singularity 
at $x=0$. $A^{MJ}(x,t)$ has an irregular singularity 
of the Poincar\'e rank three but has no other singularities. They are not 
connected by any rational transform of the independent variable.

In this paper, we show that both (MJ) and (FN) comes from 
different degeneration from the sixth Painlev\'e equations. Moreover, we will show that 
it is natural to consider (FN) is a deformation for the thirty-fourth Painlev\'e 
equation P34($\alpha$) in   Gambier's classification \cite{Gam}
$$y^{\prime\prime}=\frac{{y^{\prime }}^2}{2y}+2y^2 - ty -\frac{\alpha}{2y},$$
which is equivalent to the second Painlev\'e equation. Instead of the original P34, 
we change the sign $t \to -t$. We call this equation as P34${}^\prime$.

It is known by Garnier and Okamoto that all types of 
the Painlev\'e equations are represented as isomonodromic deformations 
of a single linear equation with order two \cite{Oka}.  (MJ) is essentially equivalent to 
the Garnier-Okamoto form. 
From the viewpoint of Garnier-Okamoto form, we obtain a well-known coalescent diagram
of the Painlev\'e equations:
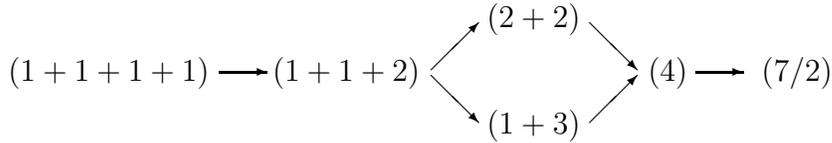
\begin{figure}[h]  
	\begin{picture}(380,55)(-40,0)
        \put(0,20){$(1+1+1+1)$}
        \put(100,20){$(1+1+2)$}
        \put(181,40){$(2+2)$}
      % \put(190,60){$(1)(1/2)$}
      % \put(265,60){$(1/2)^2$}
      % \put(120,30){$(0)^2(1/2)$}
      % \put(190,0){$(0)(3/2)$} %%
        \put(181,0){$(1+3)$}
        \put(242,20){$(4)$} %%
        \put(285,20){$(7/2)$}
        \put(80,23){\vector(1,0){18}}
      % \put(93,33){\vector(1,0){20}}
        \put(160,22){\vector(1,-1){18}}
        \put(220,42){\vector(1,-1){18}}
        \put(160,24){\vector(1,1){18}}
        \put(220,4){\vector(1,1){18}}
      % \put(168,63){\vector(1,0){20}}
        \put(260,23){\vector(1,0){18}}
      % \put(168,3){\vector(1,0){20}}
      % \put(238,63){\vector(1,0){20}}
        \end{picture}
	\label{}
\end{figure}
%% \begin{figure}[h]  
%% 	\begin{picture}(340,75)(-40,0)
%%         \put(0,30){$(0)^4$}
%%         \put(50,30){$(0)^2(1)$}
%%         \put(120,60){$(1)^2$}
%%       % \put(190,60){$(1)(1/2)$}
%%       % \put(265,60){$(1/2)^2$}
%%       % \put(120,30){$(0)^2(1/2)$}
%%       % \put(190,0){$(0)(3/2)$} %%
%%         \put(120,0){$(0)(2)$}
%%         \put(190,30){$(3)$} %%
%%         \put(265,30){$(5/2)$}
%%         \put(23,33){\vector(1,0){20}}
%%       % \put(93,33){\vector(1,0){20}}
%%         \put(93,28){\vector(1,-1){20}}
%%         \put(162,58){\vector(1,-1){20}}
%%         \put(93,38){\vector(1,1){20}}
%%         \put(162,8){\vector(1,1){20}}
%%       % \put(168,63){\vector(1,0){20}}
%%         \put(238,33){\vector(1,0){20}}
%%       % \put(168,3){\vector(1,0){20}}
%%       % \put(238,63){\vector(1,0){20}}
%%         \end{picture}
%% 	\label{}
%% \end{figure}

\par\noindent
Here $(j)$ is a pole order of the connection $A(x,t)$.  
This diagram is easy to understand and  explains coalescence of the Painlev\'e 
equations \cite{P:06}. But it seems that (FN) is out of the coalescent diagram 
since the type of singularities of (FN) is $(1+4)$.  
Later we will show a complete coalescent diagram of the Painlev\'e equations from 
 the sixth Painlev\'e equations, which contains (FN) as the type $(1+5/2)$.

Before we show the complete coalescent diagram, we will review 
the third Painlev\'e equation P3
\begin{equation}\label{0:p3}
y^{\prime\prime}= \frac{1}{ y}{y^{\prime}}^2-\frac{y^{\prime}}{ t}+
\frac{\alpha y^2 + \beta }t+\gamma y^3 + \frac{\delta }{ y}. 
\end{equation}
P3 is divided into four type 
\begin{quote}
(P3-A) $\gamma\not= 0, \delta\not=0$ \\
(P3-B) $\gamma\not= 0,  \delta=0 $ or $\gamma= 0,\ \delta\not=0$ \\
(P3-C) $\gamma=0,\ \delta=0$ \\
(P3-D) $\alpha=0,\ \gamma=0$ or $\beta=0,\ \delta=0$,
\end{quote}
Since the case (P3-D) is quadrature, we exclude the case (P3-D). 
The cases (P3-A), (P3-B) and (P3-C) are called 
the type $D_6^{(1)}$, the type $D_7^{(1)}$, the type $D_8^{(1)}$, respectively.
The meaning of type is the Dynkin diagram of the intersection form of boundary 
divisors of the Okamoto initial value spaces  \cite{S}.  
In \cite{OKSO} we show that the corresponding linear equations for $D_7^{(1)}$ and 
$D_8^{(1)}$ type has  singularities of type $(1)(1/2)$ and $(1/2)^2$. 
These three different types of the third equations are noticed by Painlev\'e \cite{P:98b}.

In the same way, the fifth Painlev\'e equation
\begin{equation}\label{0:p5}
y^{\prime\prime} = \left(\frac{1}{2y}+\frac{1}{y-1} \right){y^{\prime}}^2-\frac{1}{t}
{y^{\prime}}+\frac{(y-1)^2}{t^2}\left(\alpha y +\frac{\beta}{y}\right)
+\gamma \frac{y}{ t}+\delta \frac{y(y+1)}{y-1},
\end{equation}
has three type 
\begin{quote}
(P5-A) $\delta\not=0$ \\
(P5-B) $\delta=0, \gamma\not= 0$ \\
(P5-C) $\delta=0, \gamma=0$
\end{quote}
 (P5-A) is a generic case. In the case (P5-B),  the fifth Painlev\'e equation is equivalent 
to the third Painlev\'e equation of type $D_6^{(1)}$. In the case (P5-C),  
the fifth Painlev\'e equation is quadrature and we exclude the case (P5-C). 
We denote the case (P5-B) as deg-P5.

 % \pagebreak
\par\bigskip

In this paper, we study a complete coalescent diagram of singularity type:
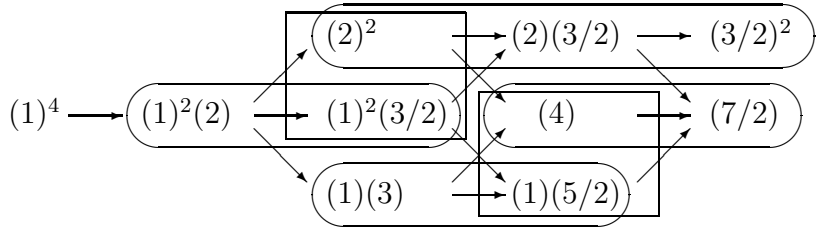
\begin{center}
\begin{figure}[h]  
	\begin{picture}(340,75)(-40,0)
        \put(0,30){$(1)^4$}
        \put(50,30){$(1)^2(2)$}
        \put(120,60){$(2)^2$}
        \put(190,60){$(2)(3/2)$}
        \put(265,60){$(3/2)^2$}
        \put(120,30){$(1)^2(3/2)$}
        \put(190,0){$(1)(5/2)$} %%
        \put(120,0){$(1)(3)$}
        \put(200,30){$(4)$} %%
        \put(265,30){$(7/2)$}
        \put(23,33){\vector(1,0){20}}
        \put(93,33){\vector(1,0){20}}
        \put(93,28){\vector(1,-1){20}}
        \put(168,58){\vector(1,-1){20}}
        \put(93,38){\vector(1,1){20}}
        \put(168,8){\vector(1,1){20}}
        \put(168,38){\vector(1,1){20}}%
        \put(168,28){\vector(1,-1){20}}%
        \put(238,8){\vector(1,1){20}}%
        \put(238,58){\vector(1,-1){20}}%
        \put(168,63){\vector(1,0){20}}
        \put(238,33){\vector(1,0){20}}
        \put(168,3){\vector(1,0){20}}
        \put(238,63){\vector(1,0){20}}
        \put(105,72){\line(1,0){68}}
        \put(105,24){\line(1,0){68}}
        \put(105,24){\line(0,1){48}}
        \put(173,24){\line(0,1){48}}
        \put(108,33){\oval(126, 24)}
        \put(175,3){\oval(120, 24)}
        \put(210,63){\oval(190 ,24)}
        \put(240,33){\oval(120, 24)}
        \put(178,42){\line(1,0){68}}
        \put(178,-5){\line(1,0){68}}
        \put(178,-5){\line(0,1){47}}
        \put(246,-5){\line(0,1){47}}
        \end{picture}
\end{figure}
\end{center}
The next diagram is the type of the Painlev\'e equation corresponding to the 
singularity diagram:
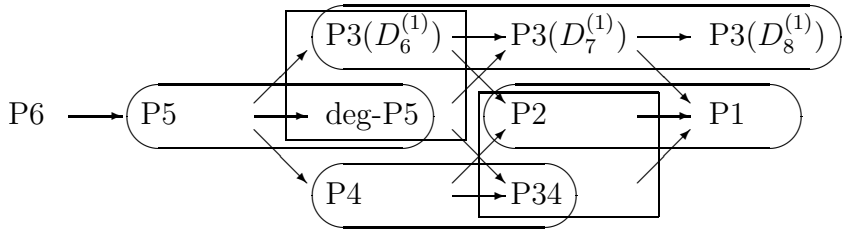
\begin{figure}[h]  
	\begin{picture}(340,75)(-40,0)
        \put(0,30){P6}
        \put(50,30){P5}
        \put(120,60){P3($D_6^{(1)}$)}
        \put(190,60){P3($D_7^{(1)}$)}
        \put(265,60){P3($D_8^{(1)}$)}
        \put(120,30){deg-P5}
        \put(190,0){P34} %%
        \put(120,0){P4}
        \put(190,30){P2} %%
        \put(265,30){P1}
        \put(23,33){\vector(1,0){20}}
        \put(93,33){\vector(1,0){20}}
        \put(93,28){\vector(1,-1){20}}
        \put(168,58){\vector(1,-1){20}}
        \put(93,38){\vector(1,1){20}}
        \put(168,8){\vector(1,1){20}}
        \put(168,63){\vector(1,0){20}}
        \put(238,33){\vector(1,0){20}}
        \put(168,3){\vector(1,0){20}}
        \put(238,63){\vector(1,0){20}}
        \put(168,38){\vector(1,1){20}}%
        \put(168,28){\vector(1,-1){20}}%
        \put(238,8){\vector(1,1){20}}%
        \put(238,58){\vector(1,-1){20}}%
        \put(105,73){\line(1,0){68}}
        \put(105,24){\line(1,0){68}}
        \put(105,24){\line(0,1){49}}
        \put(173,24){\line(0,1){49}}
        \put(103,33){\oval(116,24)}
        \put(165,3){\oval(100, 24)}
        \put(215,63){\oval(200 ,24)}
        \put(240,33){\oval(120, 24)}
        \put(178,42){\line(1,0){68}}
        \put(178,-5){\line(1,0){68}}
        \put(178,-5){\line(0,1){47}}
        \put(246,-5){\line(0,1){47}}
        \end{picture}
\end{figure}
\par\noindent
In both diagrams, we have two boxes and four ovals.
 We will show that the  Painlev\'e equations  in a box are equivalent (Theorem 1). 
The  Painlev\'e equations and  their isomonodromic deformations in an oval can be 
unified in one equation (Theorem 2). 

We add four new types  to old diagrams.  All of them have 
a singularity whose order is a half integer.  The types 
$(1)^2(3/2)$, $(1)(5/2)$, $(1)(3/2)$ and $(3/2)^2$ correspond to deg-P5, 
P34, P3($D_7^{(1)}$)   and  P3($D_8^{(1)}$), respectively. 
The third Painlev\'e equation of $D_6^{(1)}$ type and the second Painlev\'e equation have 
two different types of isomonodromic deformations. 
The type $(0)^2(1/2)$ is corresponding to the fifth Painlev\'e equation in the case $\delta=0$. 
A transformation between (FN) and type  $(1)(5/2)$ is also pointed out in \cite{KH99}.
We also add eight new arrows.  We have two types of coalescence. One is confluence of two 
singularities $(r_1)(r_2) \to  (r_1+r_2)$. The second is decrease in the Poincar\'e rank 
$(r)  \to  (r-1/2)$ when $r=2,3,4$.  In the old diagram, the second type appeared only 
in the case P2 $\to$ P1.

\begin{thm}\label{thm1}
The coalescent diagram which starts a linear differential equation with four regular singularities 
consists of ten types of singularities. We obtain eight different types of the Painlev\'e equations
from this diagram. The third Painlev\'e equation of $D_6^{(1)}$ type 
and the second Painlev\'e equation have two types of isomonodromic deformations.
\end{thm}
 
The first Painlev\'e equation P1
$$y^{\prime\prime} = 6 y^2 + t.$$
 can be considered as deg-P2.  Painlev\'e showed that a unified 
equation of P1 and P2 \cite{P:98a}:
\begin{equation}\label{P1_2}
y^{\prime\prime}= \alpha(2y^3+ty)+ \beta(6y^2+t)
%% y^{\prime\prime} = \alpha (y^3+3 t y)+ 3\beta (y^2 + t ).
\end{equation} 
In \cite{P:98a}, Painlev\'e took $\beta=1$. 
If $\alpha=0$, \eqref{P1_2} is nothing but P1. 
We will show  \eqref{P1_2} is equivalent to P2 if $\alpha\not=0$ in the section \ref{pain}.  

deg-P5 is also a special case of P5, and P3($D_7^{(1)}$) and P3($D_8^{(1)}$) are
also  special cases of P3. 
In the section \ref{pain}, we show the equation P4\_34${}^\prime$($\alpha, \beta, \gamma$) 
\begin{equation}\label{P4_34}
y^{\prime\prime}=\frac{{y^{\prime }}^2}{2y}
-\frac{\alpha}{2y} +\beta y(2y +t) +\gamma y(y+t)(3y+t).
\end{equation} 
is a unified equation of P4 and P34${}^\prime$. If $\gamma=0$, P4\_34${}^\prime$($\alpha, \beta, \gamma$) 
is  equivalent to P34${}^\prime$${}^\prime$. If $\gamma\not=0$, P4\_34${}^\prime$($\alpha, \beta, \gamma$) 
is  equivalent to P4. 
The authors cannot find the unified equation  \eqref{P4_34} in literatures.

%%  P34 can be also considered as a special case of  P4$(\alpha, \beta)$
%% $$y^{\prime\prime}= \frac{1}{2y}{y^{\prime}}^2+\frac{3}{ 2}y^3+ 4 t y^2 +2(t^2-\alpha)y +\frac{\beta}{y}.$$
%% For P4($\alpha, \beta$) we set
%% \begin{equation*} 
%%  y \to 2 \varepsilon  y,\ t\to  - \varepsilon t  + \frac {d^3}{4\varepsilon^3},\ 
%% \alpha \to   \frac {d^6}{16\varepsilon^6},\ \beta \to -2\alpha. 
%% \end{equation*} 
%% Then we obtain the unified equation P4\_34($\alpha, \beta, \gamma$) of P4 and P34
%% 
%% for $\beta=d^3, \gamma= 2\varepsilon^4$. P4\_34($\alpha, \beta, \gamma$) is equivalent to 
%% P34  if $\gamma=0$. We will study P4\_34($\alpha, \beta, \gamma$) in the section \ref{deg}
%% in detail. 
%%
%% Conversely, \eqref{P4_34} is equivalent to   P4, if $\delta\not=0$.  
%% We set $\delta= 2\varepsilon^4$.  By  changing variables
%% \begin{equation}\label{deg4to34}
%% y \to  \frac y{2 \varepsilon},\ t\to  -\frac{t}{ \varepsilon}+ \frac 1{4\varepsilon^4},\
%% \alpha \to   -\frac{\beta}4,
%% \end{equation} 
%% \eqref{P4_34} is transformed to  P4$(\varepsilon^{-6}/16, \beta)$. 
Thus we obtain the following observation.
\begin{thm}\label{thm2}
In the coalescent diagram, equations in an oval can be represented as one unified equation. 
P5 and deg-P5 are unified as the standard   fifth Painlev\'e equation.  
P($D_6^{(1)}$),   P($D_7^{(1)}$) and P($D_8^{(1)}$) are unified as the standard  third 
Painlev\'e equation.  P4 and P34 are unified as \eqref{P4_34}. 
P1 and P2 are unified as \eqref{P1_2}. 
In an oval, coalescence reduces the Poincar\'e rank of   a singularity by 1/2. 
The corresponding linear equations are also unified in one unified equation. 
\end{thm}

As a single nonlinear equations, the Painlev\'e equations are classified into five types.
Each type has a scaling transformation $t \to c_1 t, y \to c_2 y$ except P6. 
Ee obtain eight types of the Painlev\'e equations after  we classify again 
each type by the scaling transformation, 

\par\bigskip
In the section \ref{pain}, we review  the Painlev\'e equations. We will show 
that the Painlev\'e equations in the same box is equivalent. 
In the section \ref{fla}, we show that (FN) comes from an isomonodromic 
deformation of the type $(1)(5/2)$. 
We will give two types of isomonodromic deformations of the Painlev\'e equations.
One is the canonical type in the section \ref{canonical}, 
This form is easy to study when we consider ten types of the Painlev\'e equations.
And the most of the Hamiltonians are   polynomials. 
In the section \ref{deg}, we give a degeneration of the extended linear equation 
\begin{equation*} \begin{aligned}
\dfrac{d^2u}{dx^2} + & p(x,t)\dfrac{du}{dx}+q(x,t)u=0,\\
\dfrac{\partial u}{\partial t} &=a(x,t)\dfrac{\partial u}{\partial x}+b(x,t)u 
\end{aligned}\end{equation*}
of the Painlev\'e equations. 
The second is $SL$-type in the section \ref{sl}. 
In this form the extended linear equations of the Painlev\'e equations  in the same oval 
are also unified in one linear equations . But the Hamiltonians are not polynomials 
in this form.  
Most of equations and degenerations are already listed in 
\cite{Oka}, but we  correct misprints in \cite{Oka}. 

The authors thank to Professor Hiroyuki Kawamuko for fruitful discussions.

\section{List of the Painlev\'e equations}\label{pain}
In this section  we list up the Painlev\'e equations in unusual way. 
This classification is essential for our coalescent diagram. 
We also give some equivalence between different types of the 
 Painlev\'e equations. We will give a proof of the second part of the Theorem 1, 
although this is well-known.

% In this section we review some equivalence between different types of the 
% Painlev\'e equations.  We will give a proof of the second part of the Theorem 1, 
% although this is well-known.

We list {\it five} types of the Painlev\'e equations:
\begin{eqnarray*}
{\textrm P1\_2)}\quad  y^{\prime\prime}&=&\alpha (2y^3+  t y)+  \beta(6y^2 + t ),\\
{\textrm P4\_34}{}^\prime{\textrm )}\quad  y^{\prime\prime}&=& \frac{{y^{\prime }}^2}{2y}-\frac{\alpha}{2y} +\beta y(2y +  t) +\gamma y(y+  t)(3y+  t),\\
%{\textrm P4\_34)}\quad  y^{\prime\prime}&=& \frac{{y^{\prime }}^2}{2y}-\frac{\alpha}{2y} +\beta y(2y  -t) +\gamma y(y-t)(3y-t),\\
{\textrm P3)}\qquad  y^{\prime\prime}&=& \frac{1}{ y}{y^{\prime}}^2-\frac{y^{\prime}}{ t}+
\frac{\alpha y^2 +  \beta  }t+\gamma y^3 + \frac{\delta }{ y},\\
\textrm{P5)}\qquad  y^{\prime\prime}&=& \left(\frac{1}{2y}+\frac{1}{y-1} \right){y^{\prime}}^2-\frac{1}{t}
{y^{\prime}}+\frac{(y-1)^2}{t^2}\left(\alpha y +\frac{\beta}{y}\right)+\gamma {y\over t}+\delta {y(y+1)\over y-1},\\
\textrm{P6)}\qquad  y^{\prime\prime}&=& {1\over 2}\left({1\over y}+{1\over y-1}+{1 \over y-t}\right){y^{\prime}}^2
-\left({1\over t}+{1\over t-1}+{1 \over y-t}\right)y^{\prime}
\end{eqnarray*}
$$
+ {y(y-1)(y-t)\over t^2(t-1)^2}\left[\alpha+\beta {t\over y^2}+\gamma {t-1\over (y-1)^2}
+ \delta{t(t-1)\over (y-t)^2}\right].
$$
\par\bigskip\noindent
Here $\alpha ,\beta, \gamma ,\delta$ are complex parameters.
P1\_2, P4\_34${}^\prime$, P3 and P5 have a scaling transformation. We will classify five types to
fourteen types   by scaling transformations. 

\subsection{Unified equation of P1 and P2}
By the scaling transformation 
 $ y \to c y,\  t \to c^2 t,$ 
P1\_2($\alpha, \beta$)  is changed to  P1\_2($  c^6\alpha, c^5\beta$).
P1\_2($\alpha, \beta$)  is divided into three types 
\begin{quote}
(P1-A) $\alpha\not=0$, \\
(P1-B) $\alpha=0, \beta\not=0$, \\
(P1-C) $\alpha=0, \beta=0$.
\end{quote}
\begin{lemma}
The case (P1-A) is equivalent to P2 and the case (P1-B) is equivalent to P1:
\begin{eqnarray*}
{\textrm P1)}\quad  y^{\prime\prime}&=& 6y^2 + t,\\
{\textrm P2)}\quad  y^{\prime\prime}&=& 2y^3+  t y +\alpha.
\end{eqnarray*}
The case (P1-C) is trivial.
\end{lemma}
\noindent
{\it Proof.}  In the case (P1-B), we can set $\beta=1$ by a scaling transformation 
and P1\_2($0, 1$)  is nothing but P1
In the case (P1-A), we set  $\alpha = \varepsilon^6$ and 
change the variables 
$$ y\to y \varepsilon^{-1}- \beta\varepsilon^{-6},\ 
t\to  t\varepsilon^{-2}+ 6\beta^2\varepsilon^{-12}.$$
Then we obtain P2
$$y^{\prime\prime}=  2y^3+ty+ \frac{4\beta^3}{\varepsilon^{ 15}}.$$
Therefore P1\_2($\varepsilon^6, \beta$) is equivalent to P2(${4\beta^3} \varepsilon^{-15}$). 
\hfill $\boxed{}$

\subsection{Unified equation of P34 and P4}
By the scaling transformation $ y \to c y, \ t \to c  t,$ 
P4\_34${}^\prime$($\alpha, \beta, \gamma$)  is changed to  P4\_34${}^\prime$($ \alpha,$ $ c^3\beta,$ $ c^4\gamma$).
P4\_34${}^\prime$($\alpha, \beta, \gamma$)  is divided into three types 
\begin{quote}
(P4-A) $\gamma\not=0$, \\
(P4-B) $\beta\not=0, \gamma=0$, \\
(P4-C) $\beta=0, \gamma=0$.
\end{quote}
\begin{lemma}
The case (P4-A) is equivalent to P4 and the case (P4-B) is equivalent to P34:
\begin{eqnarray*}
{\textrm P34}{}^\prime{\textrm )}\quad  y^{\prime\prime}&=& \frac{{y^{\prime }}^2}{2y}+2y^2 +ty -\frac{\alpha}{2y},\\
{\textrm P4)}\quad  y^{\prime\prime}&=& \frac{1}{2y}{y^{\prime}}^2+\frac{3}{ 2}y^3+ 4 t y^2 +2(t^2-\alpha)y +\frac{\beta}{y}.
\end{eqnarray*}
P2 and P34${}^\prime$ are equivalent. The case (P4-C) is quadrature.
\end{lemma}
\noindent
{\it Proof.} 
In the  case (P4-C), P4\_34${}^\prime$($\alpha, 0, 0$) has a solution
$$y= C_1 t^2 + C_2 t+ \frac{C_2^2-\alpha}{4C_1}.$$
 In the case (P4-B), we can set $\beta=1$ by a scaling transformation 
and P4\_34${}^\prime$($\alpha, 1, 0$)  is nothing but P34${}^\prime$($\alpha$). In the case (P4-A), we set 
 $\beta = d^3,  \gamma=2\varepsilon^4$ and change the variables
$$  y\to \frac{y}{2 \varepsilon},\ 
t\to   \varepsilon^{-1} {t} - \frac{d^3}{4  \varepsilon^4},\ \alpha \to  -\beta/2.$$
Then we obtain P4($d^6\varepsilon^{-6}/16,\beta$)
$$y^{\prime\prime}= \frac{1}{2y}{y^{\prime}}^2+\frac{3}{ 2}y^3+ 4 t y^2 
+2 t^2 y-\frac{d^6 y}{8 \varepsilon^6} +\frac{\beta}{y}.$$

We will show the equivalence between P2 and P34${}^\prime$. 
The second Painlev\'e equation \eqref{0:p2} is represented by a Hamiltonian form:
\begin{equation}\label{2:ha2}
{\cal H}_{II}: \
\left\{\begin{array}{l}
   q^{\prime } = -q^2 +p -\frac t2, \\
   p^{\prime } = 2pq+a, \\
\end{array}
\right.
\end{equation}
with the Hamiltonian  
$$ H_{II}= {1 \over 2}p^2- \left(q^2+{t \over 2} \right) p-a q.$$
If we remove $p$ from \eqref{2:ha2}, we obtain P2($a-1/2)$. 
If we remove $q$ from \eqref{2:ha2}, we obtain P34($a^2$). 
Therefore P2 and P34 are equivalent. P34 and P34${}^\prime$ are equivalent by $t \to -t$.

More precisely, if 
$y$ satisfies the second Painlev\'e equation P2($\alpha$),
the function $p=y^2+y^{\prime}+ t/2$ satisfies P34($(\alpha+1/2)^2$). 
Conversely, If $p$ satisfies P34($ \alpha$), 
$q=\frac 1 {2p}\left(p^{\prime }-\sqrt{\alpha} \right)$
satisfies  P2($\sqrt{\alpha}-1/2$).

\hfill $\boxed{}$

{\it Remark.}\  If we choose  $-\sqrt{\alpha}$ instead of $\sqrt{\alpha}$, 
we obtain  P2($-\sqrt{\alpha}-1/2$) which is equivalent to P2($\sqrt{\alpha}-1/2$) 
by  a B\"acklund transformation. The equivalence of P2 and P34 are known by \cite{Gam}.

We will use P34${}^\prime$ instead of P34. If we change the sign of $t$, we 
obtain a canonical transformation $(p,q,H,t) \to (q,p,H,-t)$:
$$dp\wedge dq -dH\wedge dt = -\left(dq\wedge dp -dH\wedge d(-t) \right).$$ 
In the followings, we may use $\sigma =\pm 1$ to express the both of P34 and P34${}^\prime$: 
$$y^{\prime\prime}=\frac{{y^{\prime }}^2}{2y}+2y^2 +\sigma ty -\frac{\alpha}{2y}.$$
Similarly, we express the both of P4\_34 and P4\_34${}^\prime$:
$$y^{\prime\prime}= \frac{{y^{\prime }}^2}{2y}-\frac{\alpha}{2y} +\beta y(2y +\sigma  t) +\gamma y(y+\sigma t)(3y+\sigma t).$$

\subsection{P3}
By the scaling transformation 
$y\to c_1y,\ t\to c_2 t$, 
P3($\alpha, \beta, \gamma, \delta$) is changed to  
P3($c_1 c_2\alpha,$ $ c_2/c_1 \beta,$ $ c_1^2c_2^2 \gamma,$ $ c_2^2/c_1^2\delta$). 
P3($\alpha, \beta, \gamma, \delta$)  is divided into four types 
\begin{quote}
(P3-A) $\gamma\not= 0, \delta\not=0$ \\
(P3-B) $\gamma\not= 0,  \delta=0 $ or $\gamma= 0,\ \delta\not=0$ \\
(P3-C) $\gamma=0,\ \delta=0$ \\
(P3-D) $\alpha=0,\ \gamma=0$ or $\beta=0,\ \delta=0$.
\end{quote}
(P3-A) is P3($D_6^{(1)}$), (P3-B) is P3($D_7^{(1)}$), (P3-C) is P3($D_8^{(1)}$)
and (P3-D) is quadrature. In usual we fix $\gamma=4, \delta=-4$ for P3($D_6^{(1)}$), 
 $\alpha=2, \gamma= 0,  \delta=-4 $ for P3($D_7^{(1)}$)  and $\alpha=4, \beta=-4, \gamma= 0,  \delta=0 $
 for P3($D_8^{(1)}$). See  \cite{OKSO}.
 
 We will use another form of the third Painlev\'e equation P3${}^\prime(\alpha, \beta, \gamma, \delta)$ 
\begin{equation*}
q^{\prime\prime}=
\frac{1}{q}{q^{\prime}}^2-\frac{q^{\prime}}{x}+
\frac{\alpha  q^2+\gamma q^3} {4x^2}+\frac{\beta}{4x} + \frac{\delta }{4q},
\end{equation*}
since P3${}^\prime$ is more sympathetic to isomonodromic deformations than P3. 
We can change P3 to P3${}^\prime$ by $x=t^2, ty = q.$

\subsection{P5}
By the scaling transformation 
$  t\to c t$, 
P5($\alpha, \beta, \gamma, \delta$) is changed to  
P5($\alpha, \beta, c \gamma, c^2\delta$).
P5($\alpha, \beta, \gamma, \delta$)  is divided into three types 
\begin{quote}
(P5-A) $\delta\not=0$ \\
(P5-B) $\gamma\not= 0,  \delta=0 $\\
(P5-C) $\gamma=0,\ \delta=0$.
\end{quote}
The case (P5-A) is a generic P5 and we call  (P5-B) as deg-P5.
In usual we fix $\delta=-1/2$ for (P5-A) and $\gamma=-2,  \delta=0 $ for (P5-B). 

\begin{lemma} 
(P5-B) is equivalent to P3${}^\prime(D_6^{(1)})$  and 
(P5-C) is quadrature.
\end{lemma}
\noindent {\it Proof.}\
P3${}^\prime(D_6^{(1)})$ is represented by a Hamiltonian form:
\begin{equation}\label{2:ha3_1}
{\cal H}_{D_6}^\prime:
\left\{\begin{array}{l}
  t q^{\prime } = 2pq^2-q^2+(\alpha_1+\beta_1)q+t, \\
  t p^{\prime } = -2p^2q + 2pq-(\alpha_1+\beta_1) p +\alpha_1. \\
\end{array}
\right.
\end{equation}
with the Hamiltonian
$$t H_{D_6}^\prime =q^2 p^2-(q^2-(\alpha_1+\beta_1) q-t)p - \alpha_1 q.$$
If we eliminate $p$ from \eqref{2:ha3_1} , 
$q$ satisfies P3${}^\prime(4(\alpha_1-\beta_1),-4(\alpha_1+\beta_1-1), 4, \ -4)$.
If we eliminate $q$ from \eqref{2:ha3_1} and set   $y=1-1/p$, 
$y$ satisfies  deg-P5$( \alpha_1^2/2, -  \beta_1^2/2, -2, 0)$.
We can write down $y$ directly by $q$:
$$ y =\frac{tq'-q^2-(\alpha_1+\beta_1)q-t}{tq'+q^2-(\alpha_1+\beta_1)q-t}.$$
Therefore deg-P5 is equivalent to P3($D_6^{(1)}$). This is known by   \cite{Grom}.
\hfill $\boxed{}$

\subsection{Summary}
If we  classify the five types of the Painlev\'e equation by scaling transformations,
we obtain fourteen types of equations. Four of them are quadrature. 
Thus we have  {\it ten} types of the Painlev\'e equations:
\begin{quote}
(P1-A), (P1-B), (P4-A), (P4-B), (P3-A), (P3-B), (P3-C), (P5-A), (P5-B), (P6).
\end{quote}
(P1-A) and (P4-B) are equivalent and (P3-A) and (P5-B) are equivalent.

\section{The Flaschka-Newell form and P34}\label{fla}
In this section we will prove that (FN) comes from isomonodromic deformations 
of type $(1)(5/2)$ and show that it it natural to consider the Flaschka-Newell 
form as an isomonodromic deformation of P34 not of P2. This proves the rest part 
of the Theorem \ref{thm1}.  The relation between the Flaschka-Newell form and 
P34 are noticed by Kapaev and Hubert \cite{Kapaev}  \cite{KH99}. 

At first we will review the Poincar\'e rank of irregular singularities.
We consider a linear equation
\begin{equation}\label{2:sing}
 \frac{d^2u}{dx^2}+p_1(x)\frac{du}{dx}+p_2(x)u=0.
\end{equation}
Assume that 
$$p_1(x)=c_0x^k+c_1x^{k-1}+\cdots, \ p_2(x)=d_0 x^l+d_1 x^{l-1}+\cdots, 
$$
around $x=\infty$ and $c_0,  d_0$ are not zero. If 
$$r={\rm max}\ (k+1,(l+2)/2) $$ 
is positive,  $x=\infty$  is an irregular singularity of \eqref{2:sing}. 
We call $r$ as the Poincar\'e rank of \eqref{2:sing} at  $x=\infty$. 
The Poincar\'e rank $r$ may be a half integer.  If $x=\infty$  is an irregular singularity 
with the Poincar\'e rank $r$, \eqref{2:sing} has solutions with an asymptotics
$$u_j\sim \exp \left(\kappa_j x^r \right).$$

\begin{prop} The Flaschka-Newell form of P2 is a double cover of a linear equation 
of the singularity type $(1)(5/2)$. If we write the equation of the type $(1)(5/2)$
as a single equation, the apparent singularity  satisfies  P34. 
\end{prop}
\par\bigskip\noindent {\it Proof.}\
We consider the following deformation equation. 
\begin{equation}\label{fn:org}\begin{aligned}
\frac{dZ}{dw}&=\left[ \begin{pmatrix} 0& 2w\\ 0&0 \end{pmatrix} +
\begin{pmatrix} -2y & -y^2-z-t/2\\ 2 &2y \end{pmatrix}+
\begin{pmatrix} \ - \alpha+1/2& 0\\ -2y^2+2z-t&  \alpha-1/2 \end{pmatrix}\frac 1{2w}
\right]Z,\\
\frac{\partial Z}{\partial t}&=  \begin{pmatrix} y& -w\\ -1&-y\end{pmatrix} Z.
\end{aligned}\end{equation}
By the compatibility condition, we obtain  P2$(\alpha)$
$$y^{\prime}=z, \quad z^{\prime}= 2y^3 +t y + \alpha.
$$
If we change $w=x^2$ and $Z=RY$  where
$$R=\begin{pmatrix} \sqrt{x}& \sqrt{x}\\ -1/\sqrt{x}&1/\sqrt{x} \end{pmatrix},$$
we obtain the FN form \eqref{FN}.  
Since the exponents of \eqref{fn:org} at $w=\infty$ coincide, the Poincar\'e rank 
at $w=\infty$ in \eqref{fn:org}  is $(3/2)$. 
We will rewrite \eqref{fn:org} as a single equation of the second order.

We change the variables 
$$w \to \frac w2,\  z\to  p^2+q- \frac t2,\ y \to -p.$$
Then  \eqref{fn:org} is changed to 
\begin{equation}\label{fn:pq}\begin{aligned}
\frac{dZ}{dw}&=\left[ \begin{pmatrix} p& x/2-q/2\\ 1&-p \end{pmatrix} +
\begin{pmatrix} -\alpha/2+1/4  &0\\q-2p^2-t &\alpha/2-1/4 \end{pmatrix}\frac 1{w}
\right]Z,\\
\frac{\partial Z}{\partial t}&=  \begin{pmatrix} -p& -x/2\\ -1&p\end{pmatrix} Z.
\end{aligned}\end{equation}
The compatibility condition is
\begin{equation} \label{p34:fn}
\left\{\begin{array}{l}
  q^{\prime } =  -2pq +\left(\alpha +\dfrac12\right), \\
  p^{\prime } =  p^2-q+  \dfrac t2.
\end{array}\right.
\end{equation}

We set 
$$Z=\begin{pmatrix} u_1\\ u_2 \end{pmatrix}, \quad u_1= w^{1/4- \alpha/2}u.$$
Eliminating $u_2$ from \eqref{fn:pq}, we get  a single equation for $u=u_1$:
\begin{equation}\label{p34:sing}\begin{aligned}
\dfrac{d^2u}{dw^2} + & p_1(w,t)\dfrac{du}{dw}+p_2(w,t)u=0,\\
\dfrac{\partial u}{\partial t} &=a(w,t)\dfrac{\partial u}{\partial w}+b(w,t)u,
\end{aligned}\end{equation}
where
\begin{equation*} \begin{aligned}
p_1(w,t)=-\dfrac{1}{w-q}+\dfrac{1/2-\alpha}w, \quad & 
p_2(w,t)=-\dfrac w2 +\frac t2 +\dfrac {{\cal H}_{34}}w +\dfrac {pq}{ w(w-q)}, \\
a(w,t)=-\dfrac {w}{ w-q}, \quad & b(w,t)= \dfrac {pq}{w-q}, \\
{\cal H}_{34}=   -qp^2 + & \left(\alpha +\dfrac12 \right)p + \dfrac{q^2}2 -\dfrac12 tq. 
\end{aligned}\end{equation*}
The isomonodromic deformation is equivalent to the Hamiltonian system \eqref{p34:fn}
with the Hamiltonian ${\cal H}_{34}$. 
%% \begin{equation} \label{ham:p34}
%% \left\{\begin{array}{l}
%%   p^{\prime } = \dfrac{\partial {\cal H}_{34}}{dq} = 2pq+\left(\alpha +\dfrac12\right)q-t, \\
%%   q^{\prime }= -\dfrac{\partial {\cal H}_{34}}{dp} = -q^2 + 2p -\dfrac t2. 
%% \end{array}\right.\end{equation} 
If we eliminate $p$ from \eqref{p34:fn}, we obtain P34($(\alpha+1/2)^2$) for $q$.
\par\bigskip
The first equation of \eqref{p34:sing} has an regular singularity $w=0$ and
 an irregular singularity of the Poincar\'e rank $3/2$ at $w=\infty$. 
 It also has an apparent singularity $w=q$. When we write the Painlev\'e equations 
 as isomonodromic deformations of linear equations of the second order, they have 
 an apparent singularity. And the apparent singularity is the Painlev\'e function. 
 Moreover 
 $$p={\textrm{Res}}_{w=q}p_2(w,t)$$
 is a canonical coordinate \cite{Oka}.  In the Flaschka-Newell case, the apparent singularity 
 $q$ satisfies  P34 but not P2. 
 \hfill $\boxed{}$

\section{Isomonodromic deformations of canonical type} \label{canonical}
%% It is easy to study degeneration for the canonical type, because the most of the Hamiltonians 
%% are polynomials for the  canonical type.
This section is the revision 
of the section 4.3 in \cite{Oka}. In this part, we list up isomonodromic deformations of 
the canonical type $L_J$:
\begin{equation}\label{single:mpd}\begin{aligned} 
\frac{\partial^2 u}{\partial x^2}&+p(x,t)\frac{\partial  u}{\partial x }+q(x,t)u=0,\\
\frac{\partial  u}{\partial t}=& a(x,t)\frac{\partial  u}{\partial x }+b(x,t)u.
\end{aligned}\end{equation}
The extended linear equation $L_J$ is called the canonical type if it is obtained 
from the canonical type equation $L_{VI}$ by the process by step-by-step confluence. 
And the Fuchsian equation $L_{VI}$ is called  the canonical type if either of the local 
exponents at any singular point is zero.
The compatibility condition of \eqref{single:mpd}  is
\begin{equation*} \begin{aligned} 
p_t(x,t)-\left(p(x,t)a(x,t)\right)_x+a_{xx}(x,t)+2b_x(x,t)&=0,\\
q_t(x,t)-2q(x,t)a_x(x,t)-q_x(x,t)a(x,t)+p(x,t)b_x(x,t)&+b_{xx}(x,t)=0.
\end{aligned}\end{equation*}
The second equation is an essential deformation equation and is a Hamiltonian 
system with the Hamiltonian ${\cal H}_J$.
 $b(x,t)$ is determined from the first equation by integration and 
if we change $b(x,t) \to b(x,t)+s(t)$, the compatibility condition is also 
satisfied. We can eliminate $s(t)$ by the transformation $u \to u\exp\int s(t)dt$. 
In the following list, we may change $b(x,t)$ up to an additive term $s(t)$. 

\subsection{List  of  canonical type} 
We have {\it ten} types of isomonodromic deformations: P1, P2, P34, 
P3${}^\prime(D_6^{(1)})$, P3${}^\prime(D_7^{(1)})$, 
P3${}^\prime(D_8^{(1)})$, P4, P5, deg-P5 and P6. We will show isomonodromic deformation 
not only for P3${}^\prime$ but also for the original P3.  We need two types of 
P3${}^\prime(D_7^{(1)})$ for degeneration from deg-P5.  One is the case $\gamma=0$ and 
the other is the case $\delta=0$. 
We also show the isomonodromic deformations for P$1\_2$ and  P$4\_34$. But these unified 
equations are not  necessary for degenerations. 

We will list up seventeen types, but they are classified in ten types up to 
algebraic transformations.  We will show the degeneration diagram:
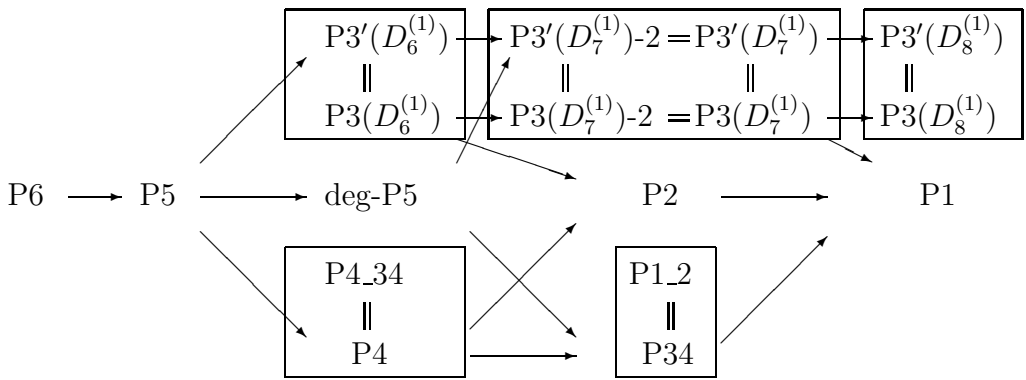
\begin{figure}[h]  
	\begin{picture}(350,135)(0,0)
        \put(0,60){P6}
        \put(50,60){P5}
        \put(120,30){P4\_34}
        \put(235,30){P1\_2}
        \put(120,90){P3($D_6^{(1)}$)}
        \put(190,90){P3($D_7^{(1)}$)-2}
        \put(260,90){P3($D_7^{(1)}$)}
        \put(330,90){P3($D_8^{(1)}$)}
        \put(120,120){P3${}^\prime(D_6^{(1)}$)}
        \put(190,120){P3${}^\prime(D_7^{(1)}$)-2}
        \put(260,120){P3${}^\prime(D_7^{(1)}$)}
        \put(330,120){P3${}^\prime(D_8^{(1)}$)}
        \put(120,60){deg-P5}
        \put(240,0){P34} %%
        \put(130,0){P4}
        \put(240,60){P2} %%
        \put(345,60){P1}
        \put(23,63){\vector(1,0){20}}
        \put(73,63){\vector(1,0){40}}
        \put(73,50){\vector(1,-1){40}}
        \put(73,76){\vector(1,1){40}}
        \put(175,50){\vector(1,-1){40}}
        \put(175,13){\vector(1,1){40}}
%%      \put(175,63){\vector(1,0){40}}
        \put(175,3){\vector(1,0){40}} 
        \put(270,8){\vector(1,1){40}}%
        \put(270,63){\vector(1,0){40}}
        \put(170,76){\vector(1,2){20}}
        \put(170,93){\vector(1,0){17}}
        \put(170,85){\vector(3,-1){45}}
        \put(170,123){\vector(1,0){17}}
        \put(250,94){\line(1,0){8}}
        \put(250,124){\line(1,0){8}}
        \put(250,92){\line(1,0){8}}
        \put(250,122){\line(1,0){8}}
        \put(135,13){\line(0,1){10}}
        \put(137,13){\line(0,1){10}}
        \put(250,13){\line(0,1){10}}
        \put(252,13){\line(0,1){10}}
        \put(135,103){\line(0,1){10}}
        \put(137,103){\line(0,1){10}}
        \put(210,103){\line(0,1){10}}
        \put(280,103){\line(0,1){10}}
        \put(340,103){\line(0,1){10}}
        \put(212,103){\line(0,1){10}}
        \put(282,103){\line(0,1){10}}
        \put(342,103){\line(0,1){10}}
        \put(310,93){\vector(1,0){17}}
        \put(310,123){\vector(1,0){17}}
        \put(310,85){\vector(2,-1){17}}
        \put(105,44){\line(1,0){68}}
        \put(105,-5){\line(1,0){68}}
        \put(105,-5){\line(0,1){49}}
        \put(173,-5){\line(0,1){49}}
        \put(105,134){\line(1,0){68}}
        \put(105,85){\line(1,0){68}}
        \put(105,85){\line(0,1){49}}
        \put(173,85){\line(0,1){49}}
        \put(182,134){\line(1,0){133}}
        \put(182,85){\line(1,0){133}}
        \put(182,85){\line(0,1){49}}
        \put(315,85){\line(0,1){49}}
        \put(324,134){\line(1,0){60}}
        \put(324,85){\line(1,0){60}}
        \put(324,85){\line(0,1){49}}
        \put(384,85){\line(0,1){49}}
        %
%
%        \put(103,33){\oval(116,24)}
%        \put(165,3){\oval(100, 24)}
%        \put(215,63){\oval(200 ,24)}
%        \put(240,33){\oval(120, 24)}
%
        \put(230,44){\line(1,0){38}}
        \put(230,-5){\line(1,0){38}}
        \put(230,-5){\line(0,1){49}}
        \put(268,-5){\line(0,1){49}}
        \end{picture}
\end{figure}
\par\noindent
Here double lines mean algebraic transformations and 
equations in a box are equivalent to each other. 

%% See the diagram in the section \ref{deg}.

\par\bigskip\noindent
Painlev\'e 1\_2:
\begin{equation*}\begin{aligned} 
p(x,t)&= -2\eta x^2 -\eta t-\frac{1}{x-y},\\
q(x,t)&=-4 \beta x^3-(\eta +2\beta t)x-2 {\cal H}_{1\_2}+\frac{z}{x-y},\\
a(x,t)&=\frac1{2(x-y)},\\
b(x,t)&= \frac{\beta\eta^{-1}-\eta y}2   - \frac{z}{2(x-y)}.
\end{aligned} \end{equation*}
$${\cal H}_{1\_2}=\frac12 z^2 -\left(\eta y^2 +\frac12 \eta t \right)z 
-2\beta y^3 -t \beta y -\frac12 \eta y.$$
$$\alpha=\eta^2.$$

\par\bigskip\noindent
Painlev\'e I:
\begin{equation*}\begin{aligned} 
p(x,t)&=-\frac{1}{x-y},\\
q(x,t)&=-4x^3-2tx-2 {\cal H}_I+\frac{z}{x-y},\\
a(x,t)&=\frac1{2(x-y)},\\
b(x,t)&=- \frac{z}{2(x-y)}.
\end{aligned} \end{equation*}
$${\cal H}_I=\frac12 z^2 -2y^3-ty.$$

\par\bigskip\noindent
Painlev\'e II: %%OK
\begin{equation*}\begin{aligned} 
p(x,t)&=-2x^2-t-\frac{1}{x-y},\\
q(x,t)&=-(2\alpha+1)x-2 {\cal H}_{II}+\frac{z}{x-y},\\
a(x,t)&=\frac1{2(x-y)},\\
b(x,t)&=-\frac{y}{2} - \frac{z}{2(x-y)}.
\end{aligned} \end{equation*}
$${\cal H}_{II}=\frac12 z^2 -\left( y^2+\frac12  t\right)z
-\left( \alpha +\frac12\right)y.$$

\par\bigskip\noindent
Painlev\'e 4\_34: If $\sigma=+1$, this gives  P4\_34${}^\prime$. If $\sigma=-1$, this gives  P4\_34.
\begin{equation*} \begin{aligned}
p(x,t)&=\theta - \sigma \eta t +\frac{1-\kappa_0}x -\eta x-\dfrac{1}{x-y}, \\
q(x,t)&= \frac{\theta^2}4+\frac{(\kappa_0-1)\eta}2 
   -\dfrac {\sigma{\cal H}_{4\_34}}x +\dfrac {yz}{x(x-y)}, \\
a(x,t)&= \dfrac {\sigma x}{ x-y}, \\
b(x,t)&=  -\dfrac {\sigma yz}{ x-y }.
\end{aligned}\end{equation*}
$${\cal H}_{4\_34}= \sigma yz^2 - \sigma \left(\eta y^2- \theta y+\kappa_0 \right)z 
  +\sigma \left(\frac{\theta^2}4+\frac{(\kappa_0-1)\eta}2 \right)y -\eta t yz.$$ 
$$\alpha=  \kappa_0^2,\ \beta=-\eta\theta,\ \gamma=\frac12\eta^2.$$

\par\bigskip\noindent
Painlev\'e 34: If $\sigma=+1$, this gives  P34${}^\prime$. If $\sigma=-1$, this gives  P34. 
\begin{equation*} \begin{aligned}
p(x,t)&=-\dfrac{1}{x-y}+\dfrac{1-\kappa_0}x, \\
q(x,t)&=-\dfrac x2 -\frac {\sigma t}2 -\dfrac {\sigma{\cal H}_{34}}x +\dfrac {yz}{x(x-y)}, \\
a(x,t)&= \dfrac {\sigma x}{ x-y}, \\
 b(x,t)&= -\dfrac {\sigma yz}{ x-y }.
\end{aligned}\end{equation*}
$${\cal H}_{34}= \sigma \left( yz^2 -\kappa_0 z -  \dfrac{y^2}2\right) -\dfrac12 ty.$$ 
$$\alpha=\kappa_0^2.$$

\par\bigskip\noindent
Painlev\'e IV:
\begin{equation*}\begin{aligned} 
p(x,t)&=\frac{1-\kappa_0}x-\frac{x+2t}2-\frac{1}{x-y},\\
q(x,t)&= \frac12 \theta_\infty-\frac{{\cal H}_{IV}}{2x}+\frac{yz}{x(x-y)},\\
a(x,t)&=\frac{2x}{ x-y },\\
b(x,t)&=-\frac1{2}\left( y+2t \right)-\frac {2yz}{x-y}.
\end{aligned} \end{equation*}
$${\cal H}_{IV}=2yz^2 -(y^2 +2ty+2\kappa_0)z+\theta_\infty  y.$$
$$\alpha= -\kappa_0+2\theta_\infty+1,\ \beta=-2\kappa_0^2.$$

\bigskip\noindent
Painlev\'e III$(D_6^{(1)})$:
\begin{equation*}\begin{aligned} 
p(x,t)&=\frac{\eta_0 t}{x^2}+\frac{1-\theta_0}x- \eta_\infty t-\frac{1}{x-y},\\
q(x,t)&=\frac{\eta_\infty(\theta_0+\theta_\infty)t}{2x}-\frac{t{\cal H}_{III}+yz}{2 x^2}+\frac{yz}{x(x-y)},\\
a(x,t)&=\frac{2yx}{t(x-y)}+\frac xt,\\
b(x,t)&=- \eta_\infty y - \frac{2y^2z}{t(x-y)}. %%
\end{aligned} \end{equation*}
$$t {\cal H}_{III}= 2y^2z^2 -\left\{2\eta_\infty ty^2 +(2\theta_0+1)y-2\eta_0 t \right\}z +\eta_\infty(\theta_0+\theta_\infty)ty.$$
$$\alpha=-4\eta_\infty \theta_\infty,\ \beta=4\eta_0(1+ \theta_0),\ \gamma=4 \eta_\infty^2,\ \delta=-4\eta_0^2.$$

\bigskip\noindent
Painlev\'e III$(D_7^{(1)})$: The case $\gamma=0$.
\begin{equation*}\begin{aligned} 
p(x,t)&=\frac{\eta_0 t}{x^2} + \frac{1-\theta_0}{x} - \frac{1}{x-y},\\
q(x,t)&=\frac {\theta_\infty t}{2x}-\frac{t{\cal H}_{D_7}+yz}{2 x^2}+\frac{yz}{x(x-y)},\\
a(x,t)&=\frac{2yx}{t(x-y)}+\frac xt,\\
b(x,t)&=-   \frac{2y^2z}{t(x-y)}. %%
\end{aligned} \end{equation*}
$$t {\cal H}_{D_7}= 2y^2z^2 -\left\{(2\theta_0+1)y-2\eta_0 t \right\}z +\theta_\infty ty.$$
$$\alpha=-4\theta_\infty,\ \beta= 4(\theta_0+1)\eta_0,\ \gamma=0,\ \delta=-4\eta_0^2.$$

\bigskip\noindent
Painlev\'e III$(D_7^{(1)})$-2: The case $\delta=0$.
\begin{equation*}\begin{aligned} 
p(x,t)&=-\eta_\infty t +\frac1x - \frac{1}{x-y},\\
q(x,t)&= \frac {\theta_0 t}  {2x^3}-\frac{t{\cal H}_{D_7{\textrm -}2}+yz}{2 x^2} +\frac{\theta_\infty \eta_\infty t }{2x}+ \frac{yz}{x(x-y)},\\
a(x,t)&=\frac{2yx}{t(x-y)}+\frac xt,\\
b(x,t)&=-   \frac{2y^2z}{t(x-y)}. %%
\end{aligned} \end{equation*}
$$t {\cal H}_{D_7{\textrm -}2}= 2y^2z^2 -\left\{2 \eta_\infty t y^2+y \right\}z+ \theta_\infty \eta_\infty t y+ \frac{\theta_0 t}y.$$
$$\alpha=-4 \theta_\infty\eta_\infty,\ \beta= 4 \theta_0,\ \gamma= 4\eta_\infty^2,\ \delta=0.$$

\noindent
\par\noindent
Painlev\'e III${}(D_8^{(1)})$:
\begin{eqnarray*}
 && p(x,t) = \frac{2}{x} - \frac{1}{x-y}, \\
 && q(x,t) =   \frac{t}{2 x^3} - \frac{t {\cal H}_{D_8} +yz}{ 2x^2}
                      + \frac{t}{2 x}+ \frac{yz}{x(x-y)},\\
 && a(x,t)=\frac{2yx}{t(x-y)}+\frac xt,\\ 
 && b(x,t)=  -\frac{2y^2z }{t(x-y)},
\end{eqnarray*}
$$t {\cal H}_{D_8}= 2y^2z^2 +yz +ty+\frac{t}{y}.$$
$$\alpha=-4,\ \beta= 4,\ \gamma=0,\ \delta=0.$$

\noindent
Painlev\'e III${}^\prime(D_6^{(1)})$:
\begin{equation*}\begin{aligned} 
p(x,t)&=\frac{\eta_0 t}{x^2}+\frac{1-\theta_0}x- \eta_\infty  -\frac{1}{x-y},\\
q(x,t)&=\frac{\eta_\infty(\theta_0+\theta_\infty) }{2x}-\frac{t{\cal H}_{D_6}^\prime}{x^2}+\frac{yz}{x(x-y)},\\
a(x,t)&=\frac{ yx}{t(x-y)},\\
b(x,t)&=- \frac{y^2z}{t(x-y)}. %%
\end{aligned} \end{equation*}
$$t {\cal H}_{D_6}^\prime=  y^2z^2 -\left\{ \eta_\infty y^2 + \theta_0 y- \eta_0 t \right\}z 
+\frac12 \eta_\infty(\theta_0+\theta_\infty) y.$$
$$\alpha=-4\eta_\infty \theta_\infty,\ \beta=4\eta_0(1+ \theta_0),\ \gamma=4 \eta_\infty^2,\ \delta=-4\eta_0^2.$$

\par\bigskip\noindent
Painlev\'e III${}^\prime(D_7^{(1)}$): The case $\gamma=0$
\begin{eqnarray*}
 && p(x,t) = \frac{\eta_0 t}{x^2} + \frac{1-\theta_0}{x} - \frac{1}{x-y}, \\
 && q(x,t) = - \frac {t  {\cal H}_{D_7}^{\prime}}{x^2}+\frac{\theta_\infty}{2x}+ \frac{yz}{x(x-y)}, \\
 && a(x,t)=\frac{ yx}{t(x-y)},\\ 
 && b(x,t)=-\frac{ y^2z }{t(x-y)},
\end{eqnarray*}.
$$ t {\cal H}_{D_7}^{\prime}= y^2z^2 +(-\theta_0 y+\eta_0 t)z +\frac {\theta_\infty}2 y.$$
$$\alpha=-4\theta_\infty,\ \beta= 4(\theta_0+1)\eta_0,\ \gamma=0,\ \delta=-4\eta_0^2.$$

\par\bigskip\noindent
Painlev\'e III${}^\prime(D_7^{(1)}$)-2: The case $\delta=0$
\begin{eqnarray*}
 && p(x,t) = -\eta_\infty   + \frac{1}{x} - \frac{1}{x-y}, \\
 && q(x,t) = \frac{\theta_0 t}{2 x^3}+\frac{\theta_\infty\eta_\infty}{2x}
  - \frac {t  {\cal H}_{D_7{\textrm -}2}^{\prime}}{x^2}+ \frac{yz}{x(x-y)}, \\
 && a(x,t)=\frac{ yx}{t(x-y)},\\ 
 && b(x,t)=-\frac{ y^2z }{t(x-y)},
\end{eqnarray*}.
$$ t {\cal H}_{D_7{\textrm -}2}^{\prime}= y^2z^2 
-\eta_\infty y^2z +\frac {\theta_\infty\eta_\infty  }2 y+\frac{\theta_0 t}{2y}.$$
$$\alpha=-4 \theta_\infty\eta_\infty ,\ \beta= 4 \theta_0,\ \gamma= 4\eta_\infty^2,\ \delta=0.$$

\par\bigskip\noindent
Painlev\'e III${}^\prime(D_8^{(1)})$:
\begin{eqnarray*}
 && p(x,t) = \frac{2}{x} - \frac{1}{x-y}, \\
 && q(x,t) =  \frac{t}{2 x^3} - \frac{t {\cal H}_{D_8}^{\prime} }{ x^2}
                      + \frac{1}{2 x}+ \frac{yz}{x(x-y)},\\
 && a(x,t)=\frac{ yx}{t(x-y)},\\ 
 && b(x,t)=  -\frac{ y^2z }{t(x-y)},
\end{eqnarray*}
$$t {\cal H}_{D_8}^{\prime}= y^2z^2 +yz +\frac y2 + \frac{t}{2y}.$$
$$\alpha=-4,\ \beta=4,\ \gamma=0,\ \delta=0.$$

\bigskip\noindent
 Painlev\'e V:
\begin{equation*}\begin{aligned} 
p(x,t)&=\frac{1-\kappa_0}x+\frac{\eta t}{(x-1)^2}+\frac{1-\theta}{ x-1 }-\frac{1}{x-y},\\
q(x,t)&= \frac\kappa{x(x-1)}-\frac{t{\cal H}_{V}}{x(x-1)^2}+\frac{y(y-1)z}{x(x-1)(x-y)},\\
a(x,t)&=\frac {y-1} t \frac{x(x-1)}{ x-y },\\
b(x,t)&=  %%  -\frac{(\kappa_0+\theta)(y-1)}{2t}
-  \frac{y(y-1)^2z}{t(x-y)}.
\end{aligned} \end{equation*}
$$t {\cal H}_{V}= y(y-1)^2z^2 -\left\{\kappa_0(y-1)^2+\theta y(y-1)-\eta ty\right\}z +\kappa (y-1).$$
$$\alpha=\frac12 \kappa_\infty^2,\ \beta=-\frac12 \kappa_0^2,\ 
\gamma=-(1+\theta)\eta,\ \delta=-\frac12 \eta^2,$$
$$\kappa = \frac14 (\kappa_0 +\theta )^2-\frac14 \kappa_\infty^2.$$

\bigskip\noindent
 Painlev\'e deg-V:
\begin{equation*}\begin{aligned} 
p(x,t)&=\frac{1}{x-1}+ \frac{1-\kappa_0}x  -\frac{1}{x-y},\\
q(x,t)&=\frac{\gamma t}{2(x-1)^3}+\frac{\kappa}{x(x-1)}
-\frac{t{\cal H}_{Vd}}{x(x-1)^2}+\frac{ y(y-1)z}{ x(x-1)(x-y)},\\
a(x,t)&=\frac {y-1} t \frac{x(x-1)}{ x-y },\\
b(x,t)&= %% -\frac{ \kappa_0(y-1) }{2t}
-  \frac{y(y-1)^2z}{t(x-y)}.
\end{aligned} \end{equation*}
$$t {\cal H}_{Vd}= y(y-1)^2z^2 - \kappa_0(y-1)^2z +\kappa (y-1)+
\frac{\gamma ty}{2(y-1)}.$$
$$\alpha=\frac12 \kappa_\infty^2,\ \beta=-\frac12 \kappa_0^2,\ 
 \delta= 0,\quad  \kappa =-\frac12(\alpha+\beta)= \frac14 (\kappa_0^2 - \kappa_\infty^2 ).$$

\bigskip\noindent
Painlev\'e VI: 
\begin{equation*}\begin{aligned} 
p(x,t)&=\frac{1-\kappa_0}x+\frac{1-\kappa_1}{x-1}+\frac{1-\theta}{ x-t}-\frac{1}{x-y},\\
q(x,t)&= \frac\kappa{x(x-1)}-\frac{t(t-1){\cal H}_{VI}}{x(x-1)(x-t)}+\frac{y(y-1)z}{x(x-1)(x-y)},\\
a(x,t)&=\frac{y-t}{t(t-1)}\frac{x(x-1)}{ x-y },\\
b(x,t)&= %% \frac{(1-\kappa_0-\kappa_1-\theta)(y-t)}{2t(t-1)}
-  \frac{y(y-1)(y-t)z}{t(t-1)(x-y)}.
\end{aligned} \end{equation*}
\begin{equation*}\begin{aligned}
t(t-1)& {  {\cal H}}_{VI}= y(y-1)(y-t) z^2 \\
-&\left\{\kappa_0(y-1)(y-t)+\kappa_1 y(y-t)+(\theta -1)y(y-1)\right\}z +\kappa (y-t).
\end{aligned} \end{equation*}
$$\alpha=\frac12 \kappa_\infty^2,\ \beta=-\frac12 \kappa_0^2,\ 
\gamma=\frac12 \kappa_1^2,\ \delta=\frac12 (1-\theta)^2,$$
$$\kappa = \frac14 (\kappa_0+\kappa_1+\theta-1)^2-\frac14 \kappa_\infty^2.$$

\subsection{Degeneration}\label{deg}
In this section we   list up all of degeneration of the Painlev\'e equations and 
extended linear equations $L_J$.  Here we consider degeneration of 
the extended linear system, which includes a deformation equation.  In some cases, 
we should take a change of the dependent variable $u \to f(x, t, \varepsilon)u$. 
In \cite{Oka} Okamoto did not treat the extended linear equations. If
$b(x,t)$ is changed up to a function $r(t)$ in the limit, we denote
$b \to b+ r$. 

\par\bigskip\noindent
P6$\to$P5: %%OK
We change the variables
$$t \to 1+\varepsilon t,\ 
\kappa_1 \to \varepsilon^{-1}\eta  + \theta + 1,\ \theta \to - \varepsilon^{-1}\eta,$$
$$(\alpha \to \alpha, \beta\to \beta,\ \gamma \to -\delta \varepsilon^{-2}+\gamma \varepsilon^{-1},\ 
\delta\to \delta\varepsilon^{-2}).$$
In the limit  $\varepsilon \to 0 $, $L_{VI}$ goes to $L_{V}$ and
$${\cal H}_{VI} \to  \varepsilon^{-1} {\cal H}_{V} +O(\varepsilon^0), \qquad (\varepsilon \to 0).$$

\par\bigskip\noindent
P5$\to$deg-P5:   
We change the variables
$$ 
z \to z+ \frac{\gamma}{2\varepsilon(y-1)},\ 
\eta \to \varepsilon, \ 
\theta \to \gamma \varepsilon^{-1}.$$
Then 
$$ {\cal H}_{V} +\frac{\theta^2}{4t}
\to {\cal H}_{Vd} +O(\varepsilon^{1}), \qquad (\varepsilon \to 0).$$
For $L_{V}$  we change
$$ u \to (x-1)^{\theta/2}  u $$
and $b \to b-\theta (y-1)/(2t)$ at first.  Then  $L_{V}$ goes to $L_{Vd}$ in the limit  $\varepsilon \to 0 $.
 
\par\bigskip\noindent
P5$\to$P4:  %%OK
We change the variables
$$
t \to 1+\sqrt{2}\varepsilon t,\  
y\to \frac\varepsilon{\sqrt{2}} y,\ 
z\to \sqrt{2}\varepsilon^{-1} z,\ x \to \frac\varepsilon{\sqrt{2}} x, $$ 
$$
\kappa_\infty \to \varepsilon^{-2},\ 
\theta \to \varepsilon^{-2} + 2 \theta_\infty -\kappa_0,\ 
\eta \to -\varepsilon^{-2},\ 
$$
$$(\alpha \to  \varepsilon^{-4}/2,\ \beta\to \beta/4,\ 
 \gamma\to -  \varepsilon^{-4},\ \delta\to - \varepsilon^{-4}/2+\alpha \varepsilon^{-2}),$$
After changing variables, we set $u \to \exp (\varepsilon^{-1}t/\sqrt{2}) u$. Then 
in the limit  $\varepsilon \to 0 $, $L_{V}$ goes to $L_{IV}$ with $b\to b+t+y/2$ and
$$\sqrt{2} \left( {\cal H}_{V} +\frac{(\kappa_0+\theta)^2-\kappa_\infty^2}4 \right)
\to \varepsilon^{-1} \left( {\cal H}_{IV}+2\theta_\infty t \right)+O(\varepsilon^{0}), \qquad (\varepsilon \to 0).$$
 
\par\bigskip\noindent
P5$\to$P3${}^\prime(D_6^{(1)})$: %%OK
We change the variables
$$y\to 1+\varepsilon y, \ z \to z/\varepsilon,\ x\to 1+\varepsilon x,$$
$$
\kappa_0 \to \varepsilon^{-1} \eta_\infty,\ 
\kappa_\infty \to \varepsilon^{-1}\eta_\infty-\theta_\infty, \ 
\theta \to \theta_0,\ 
\eta \to  \varepsilon \eta_0,\ 
$$
$$
(\alpha \to \frac18{\varepsilon^{-2}\gamma}+\frac14\varepsilon^{-1}\alpha , 
\beta\to -\frac{\varepsilon^{-2}\gamma}8, \gamma\to   \frac{\varepsilon\beta}4, 
\delta\to \frac{\varepsilon^{2}\delta}8).$$
In the limit  $\varepsilon \to 0 $, $L_{V}$ goes to $L_{D_6}$ and
$$ {\cal H}_{V}  \to  {\cal H}_{D_6}^\prime +O(\varepsilon^{1}), \qquad (\varepsilon \to 0).$$

\par\bigskip\noindent
deg-P5$\to$P3${}^\prime(D_7^{(1)})$-2:
We change the variables
$$ 
y\to 1+\varepsilon  y, \ z \to \varepsilon^{-1}z,\  x\to 1+\varepsilon  x,$$
$$
\kappa_\infty \to    \varepsilon^{-1} \eta_\infty,\ 
\kappa_0 \to    \varepsilon^{-1} \eta_\infty+\theta_\infty,\ 
\gamma \to \theta_0 \varepsilon/4, \ 
$$
$$
(\alpha \to -\varepsilon^{-2}\gamma/8,\ \beta \to  \varepsilon^{-1} \alpha/4,\ 
\gamma \to \beta \varepsilon/4).$$
In the limit  $\varepsilon \to 0 $, $L_{Vd}$ goes to $L_{D_7{\textrm -}2}$ 
  %% with $b\to b -\eta_\infty y/(2t)$ 
and 
$$ {\cal H}_{Vd}  \to  {\cal H}_{D_7{\textrm -}2} +O(\varepsilon), \qquad (\varepsilon \to 0).$$

\par\bigskip\noindent
deg-P5$\to$P34: 
We change the variables
$$t \to 1+\sigma \varepsilon^2 t,\ 
y\to \varepsilon^2 y, \ z \to z\varepsilon^{-2},\  x\to \varepsilon^2 x,$$
$$
\kappa_\infty \to \sigma \sqrt{-2}\varepsilon^{-3 },\ 
\gamma \to \varepsilon^{-6}, \ 
$$
$$
(\alpha \to -\varepsilon^{-6},\ \beta \to -\alpha/2,\  \gamma \to \varepsilon^{-6}).$$
In the limit  $\varepsilon \to 0 $, $L_{Vd}$ goes to $L_{34}$  
and
$$ {\cal H}_{Vd}   =   \varepsilon^{-2}\left(  \sigma{\cal H}_{34}-t^2/2 \right)
                    +\sigma t \varepsilon^{-4}  /2 -\varepsilon^{-6}/2  +O(\varepsilon^{-1}), \qquad (\varepsilon \to 0).$$

\par\bigskip\noindent
P3${}^\prime \to$P3: %%OK
This transformation is algebraic and we do not take any limit.
If we change the variables
$$t \to t^2, y \to ty,\ z  \to z/t, \quad x \to tx,$$
$L_{J}^\prime$ is changed to $L_{J}$ with $b\to b +\eta_\infty y$ if $J=D_6^{(1)}$ and
$${\cal H}_{J}^\prime    =   
\frac1{2t} {\cal H}_{J} + \frac{yz}{2t^2},$$
for $J=D_6^{(1)}, D_7^{(1)}, D_7^{(1)}{\textrm -2}, D_8^{(1)}$.

\par\bigskip\noindent
P3${}^\prime(D_6^{(1)})  \to$P3${}^\prime(D_7^{(1)})$:
We change the parameters 
$$\eta_\infty \to \varepsilon,\  
  \theta_\infty \to  \theta_\infty\varepsilon^{-1}.$$
In the limit  $\varepsilon \to 0 $, $L_{D_6}^\prime$ goes to $L_{D_7}^\prime$ and 
$${\cal H}_{D_6}^\prime  =  {\cal H}_{D_7}^\prime +O(\varepsilon), \qquad (\varepsilon \to 0).$$
P3${}(D_6^{(1)})  \to$P3$(D_7^{(1)})$ is as the same.

\par\bigskip\noindent
P3${}^\prime(D_6^{(1)})  \to$P3${}^\prime(D_7^{(1)})$-2:
We change the parameters 
$$\eta_0 \to \varepsilon,\ \theta_0 \to \theta_0\varepsilon^{-1},\  
z \to z+ \frac{\theta_0}{2  \varepsilon y}.$$
Then we have 
$${\cal H}_{D_6}^\prime +\frac{\theta_0^2}{4 t } \to  {\cal H}_{D_7{\textrm -2}}^\prime +O(\varepsilon), \qquad (\varepsilon \to 0).$$
For $L_{D_6}^\prime$, we change 
$$u \to x^{\theta_0/2} u$$
at first. Then $L_{D_6}^\prime$  goes to $L_{D_7{\textrm -2}}^\prime$. 
P3${}(D_6^{(1)})  \to$P3$(D_7^{(1)})$-2 is as the same.

\par\bigskip\noindent
P3${}^\prime(D_7^{(1)})  \to$P3${}^\prime(D_7^{(1)})$-2:
This transformation is algebraic and we do not take any limit.
If we change the variables
$$  y \to t/y,\ z \to (\theta_\infty y/2 -y^2 z)/t,$$
$$ \theta_0 \to \theta_\infty-1,\ \theta_\infty\to  \theta_0,\ \eta_0 \to \eta_\infty, $$
 $$(\alpha \to -\beta,\ \beta \to -\alpha,\ \delta \to -\gamma). $$ 
we have
$${\cal H}_{D_7}^\prime = {\cal H}_{D_7{\textrm -2}}^\prime -\frac{yz}t 
-\frac{\theta_\infty(\theta_\infty-2)}{4t}.$$
For $L_{D_7}^\prime$ we 
change 
$$  u \to u x^{(\theta_0+1)/2} $$ %% \exp \int \frac{2yz- \theta_\infty}{2t} dt
at first. Changing the variable   $x \to     t/x,$ 
$L_{D_7}$ is changed to $L_{D_7{\textrm -2}}^\prime$ with $b \to b-yz/t$.

\par\bigskip\noindent
P3${}^\prime(D_7^{(1)})  \to$P3${}^\prime(D_8^{(1)})$:
We change the variables
$$t\to 2t,\ y \to 2y,\ z \to z/2 +1/(4\varepsilon  y),$$
$$ \eta_0 \to \varepsilon,\     \theta_0 \to -1+\varepsilon^{-1},\ \theta_\infty \to 1/2.$$
Then 
$${\cal H}_{D_7}^\prime  = \frac12 {\cal H}_{D_8}^\prime -\frac{1}{8\varepsilon^2 t}
+\frac{1}{4 \varepsilon t}  +O(\varepsilon), \qquad (\varepsilon \to 0).$$
For $L_{D_7}^\prime$  we change
$$ u \to x^{(1+\theta_0)/2}u $$
at first. Changing the variable  $x\to 2x,$ $L_{D_7}^\prime$ goes to $L_{D_8}^\prime$ 
in the limit  $\varepsilon \to 0 $.
P3${}(D_7^{(1)})  \to$P3$(D_8^{(1)})$ is as the same.

\par\bigskip\noindent
P3${}^\prime(D_7^{(1)})$-2  $\to$P3${}^\prime(D_8^{(1)})$:
We change the variables
$$t\to -2t,\  \ z \to z  +1/(2y),$$
$$ \eta_\infty \to \varepsilon,\  \theta_0 \to -1/2,\ \theta_\infty \to \varepsilon^{-1}.$$
Then 
$${\cal H}_{D_7{\textrm -2}}^\prime  = -\frac12 {\cal H}_{D_8}^\prime -\frac{1}{8  t}
+O(\varepsilon), \qquad (\varepsilon \to 0).$$
For $L_{D_7{\textrm -2}}^\prime$  we change
$$ u \to x^{1/2}u $$
at first. Changing the variables,  $L_{D_7{\textrm -2}}^\prime$ goes to $L_{D_8}^\prime$ 
in the limit  $\varepsilon \to 0 $.
P3${}(D_7^{(1)})$-2  $\to$P3$(D_8^{(1)})$ is as the same.

\par\bigskip\noindent
P3$(D_6^{(1)})\to$P2: %%OK
We change the variables
$$t \to 1+\varepsilon^2 t,\ y \to 1+2\varepsilon y,\ z\to 1+\frac{z}{2\varepsilon},\ x \to 1+2\varepsilon x, $$
$$\eta_0 \to -\varepsilon^{-3} /4,\ \eta_\infty \to \varepsilon^{-3} /4,\ 
\theta_0 \to -\varepsilon^{-3} /2 -2\alpha-1,\ \theta_\infty \to \varepsilon^{-3}/2,$$
$$(\alpha \to -\frac{\varepsilon^{-6}}2, \beta\to \frac12  \varepsilon^{-6}(1+4\alpha  \varepsilon^3), 
\gamma\to \frac{\varepsilon^{-6}}4, \delta\to -\frac{\varepsilon^{-6}}4).$$
I, 
After changing variables, we set $ \ u \to \exp(-\varepsilon^{-1}t/4)u.$ Then 
in the limit  $\varepsilon \to 0 $,  $L_{D_6}$ goes to $L_{II}$ and
$${\cal H}_{D_6} +\eta_0(\theta_0+\theta_\infty) \to   
 \varepsilon^{-2} {\cal H}_{II} +O(\varepsilon^{-1}), \qquad (\varepsilon \to 0).$$

\par\bigskip\noindent
P4$\to$P4\_34: This transformation is algebraic and we do not take any limit.
We change the variables
$$
t \to \sigma \varepsilon t - \varepsilon^{-1}\theta/2, \
y \to 2 \varepsilon y,\ 
z \to \varepsilon^{-1}z/2,\ $$
$$x \to 2 \varepsilon x,\ \theta_\infty \to (\kappa_0-1)/2 +\varepsilon^{-2}\theta^2/8,$$
$$\left(\alpha \to  \varepsilon^{-6}\beta^2/16,\ \beta \to -2\alpha,\ 2\varepsilon^4 \to \gamma \right).$$
Then    $L_{IV}$ is changed  to $L_{4\_34}$ with  
$b\to b-\sigma \varepsilon^2y-\varepsilon^2 t+\sigma \theta/2$ and
$$ {\cal H}_{IV}  =  \sigma  \varepsilon^{-1} {\cal H}_{4\_34}$$
for $\eta =2 \varepsilon^2$.

\par\bigskip\noindent
P4$\to$P2: %%OK
We change the variables
$$
t \to -\varepsilon^{-3} (1-2^{-2/3}\varepsilon^4 t), \
y \to \varepsilon^{-3} (1+2^{ 2/3}\varepsilon^2 y),\ 
z \to 2^{-2/3}\varepsilon z,\ $$
$$x \to \varepsilon^{-3} (1+2^{ 2/3}\varepsilon^2 x),\ 
u \to \exp(2^{-5/3}\varepsilon^{-2} t)u, $$
$$\kappa_0 \to \varepsilon^{-6}/2, \theta_\infty \to -\alpha-1/2, \quad
\left(\alpha \to -2\alpha- \frac1{2\varepsilon^{ 6}},\ \beta\to -\frac1{2\varepsilon^{12}}\right).$$
In the limit  $\varepsilon \to 0 $, $L_{IV}$ goes to $L_{II}$ and
$${\cal H}_{IV} \to  2^{ 2/3} \varepsilon^{-1} {\cal H}_{II} 
-\varepsilon^{-3}(\alpha+1/2)+O(\varepsilon^0), \qquad (\varepsilon \to 0).$$

\par\bigskip\noindent
P4$\to$P34:  
We change the variables
$$
t \to \varepsilon t+ \sigma \varepsilon^{-3} /4, \
y \to 2 \sigma \varepsilon y,\ 
z \to \sigma \varepsilon^{-1} z/2+\sigma \varepsilon^{-3}/8,\ $$ 
$$x \to \sigma \varepsilon x, u \to \exp\frac{\sigma x}{8\varepsilon^3}u,$$
$$ \theta_\infty \to  \varepsilon^{-6}/32, \quad
 \left(\alpha \to \varepsilon^{-6}/16, \   \beta \to  - 2\alpha  \right).$$
In the limit  $\varepsilon \to 0 $, $L_{IV}$ goes to $L_{34}$ and
$${\cal H}_{IV}   \to  
 \frac1{ \varepsilon } {\cal H}_{34} -\frac{\sigma \kappa_0}{4\varepsilon^3} +O(\varepsilon^{0}), \qquad (\varepsilon \to 0).
$$

\par\bigskip\noindent
P3$(D_7^{(1)})\to$P1:
We change the variables
$$
t \to (\varepsilon^{-10} + \varepsilon^{-6}  t)/2, \
y \to  1+  2 \varepsilon^2 y,\ 
z \to  \varepsilon^{-2} z/2-\frac{\varepsilon^{-5}}{4\sqrt{-2}} -\frac{\varepsilon^{-3}y}{2\sqrt{-2}},\ $$ 
$$ \theta_0 \to -1 +3\sqrt{-2}\varepsilon^{-5}/4,\ \theta_\infty \to -1/2,\ \eta_0\to \sqrt{-2}\varepsilon^{5},$$
$$  \left(\alpha \to 2, \
  \beta \to -6,\ \delta \to 8\varepsilon^{10} \  \right).$$
Then  
$${\cal H}_{D_7}-\frac 58 \eta_0^2 t-\frac{\eta_0}4-\frac14  \to  
 {2\varepsilon^{6}} {\cal H}_{I}  +O(\varepsilon^{7}), \qquad (\varepsilon \to 0).
$$
For $L_{D_7}$, we change
$$  u \to x^{(\theta_0-1)2} \exp \left(\frac{ \eta_0 t}{2x}-\frac{3\eta_0t}4\right) u$$
at first. Changing the variable
$$x \to 1+  2 \varepsilon^2 x,$$
$L_{D_7}$ goes to $L_I$ in the limit $\varepsilon \to 0$.

\par\bigskip\noindent
P3$(D_7^{(1)})$-2 $\to$P1:
We change the variables
$$
t \to (\varepsilon^{-10} + \varepsilon^{-6}  t)/2, \
y \to  1-  2 \varepsilon^2 y,\ 
z \to - \varepsilon^{-2} z/2-\frac{\varepsilon^{-5}}{2\sqrt{-2}},\ $$ 
$$ \theta_0 \to -1/2,\ \theta_\infty \to  1/2 -3 \varepsilon^{-5}/(2\sqrt{-2}),\ 
   \eta_\infty\to \sqrt{-2}\varepsilon^{5},$$
$$  \left(
  \alpha \to -6,\ \beta\to -2,\ \gamma \to -8\varepsilon^{10} \  \right).$$
Then  
$${\cal H}_{D_7{\textrm -2}}+\frac{\eta_\infty^2 t}2  \to  
 {2\varepsilon^{6}} \left({\cal H}_{I}  \right) -2 +O(\varepsilon^{7}), \qquad (\varepsilon \to 0).
$$
For $L_{D_7{\textrm -2}}$, we change
$$  u \to x^{(\theta_0-1)2} \exp \left(\frac{ 2 \theta_\infty x }{3}-\frac{3\eta_0t}2 \right) u$$
at first. Changing the variable
$$x \to 1- 2 \varepsilon^2 x,$$
$L_{D_7{\textrm -2}}$ goes to $L_I$ in the limit $\varepsilon \to 0$.

\par\bigskip\noindent
P34$\to$P1:
We change the variables
$$
t \to -\sigma \varepsilon^{2} t + 6\sigma \varepsilon^{-10}, \
y \to 2\varepsilon^{-4} y -2 \varepsilon^{-10},\ 
z \to  \varepsilon^{4} z/2 +\varepsilon y +\varepsilon^{-5},\ $$ 
$$\kappa_0 \to -4{\varepsilon^{-15}} \qquad \left(   \alpha \to 16{\varepsilon^{-30}}\right).  $$
Then  
$${\cal H}_{34}  \to  -\sigma\varepsilon^{-2} {\cal H}_{I} 
+6\sigma\varepsilon^{-20}-\sigma\varepsilon^{-8} t+O(\varepsilon^{-1}), \qquad (\varepsilon \to 0).
$$
For $L_{34}$, we change
$$u \to x^{ \kappa_0/2}u$$
at first. Changing the variable
$$x \to 2\varepsilon^{-4} x -2 \varepsilon^{-10}.$$
Then $L_{34}$ goes to $L_I$ in the limit $\varepsilon \to 0$.

\par\bigskip\noindent
P1\_2$\to$P2:  
This transformation is algebraic and we do not take any limit.
We change the variables
$$
t \to \varepsilon^{-2} t +6  \varepsilon^{-2}\theta^2, \
y \to \varepsilon^{-1} y -  \varepsilon^{-1}\theta,\ 
z \to \varepsilon z -2\varepsilon  \theta y +4 \varepsilon  \theta^2,\ $$
$$\eta= \varepsilon^3, \beta=\varepsilon^5\theta, \qquad \left(\alpha =  \varepsilon^6,\  \beta  =  \varepsilon^5\theta  \right).$$
Then $$ {\cal H}_{1\_2}  =   \varepsilon^{2}\left( {\cal H}_{II} 
+\frac{\theta}2- \theta^2 t \right) $$
for $\alpha=4\theta^3$ in P2.
For    $L_{1\_2}$, we change
$$u \to \exp\left(-\frac{\beta x^2}\eta+\frac{2\beta^2 x}{\eta^3}\right)u$$
at first. Changing the variable 
$$x \to \varepsilon^{-1} x -  \varepsilon^{-1}\theta,$$
$L_{1\_2}$ is changed  to $L_{II}$ for $\alpha =4\theta^3$.

\par\bigskip\noindent
P2$\to$P1:
We change the variables
$$
t \to \varepsilon^{2} t -6\varepsilon^{-10}, \
y \to \varepsilon y +\varepsilon^{-5},\ 
z \to  \varepsilon^{-1} z +(\varepsilon y +\varepsilon^{-5})^2 +(\varepsilon^{2} t -6\varepsilon^{-10})/2,\ $$ 
$$\alpha \to 4{\varepsilon^{-15}}.  $$
Then  
$${\cal H}_{II}  = \varepsilon^{-2} {\cal H}_{I} 
-6\varepsilon^{-20}+\varepsilon^{-8} t-\varepsilon^{-5}/2+O(\varepsilon), \qquad (\varepsilon \to 0).
$$
For $L_{II}$, we change
$$u \to \exp \left(\frac{x^3}3 +\frac {tx}2\right)u$$
at first. Changing the variable
$$x \to \varepsilon x +\varepsilon^{-5},$$
$L_{II}$ goes to $L_I$ in the limit $\varepsilon \to 0$.

\par\bigskip

\par\bigskip\noindent
{\it Remark.}\ In \cite{Oka}, there is a misprint in P3 $\to$P2.

\section{Isomonodromic deformations of $SL$-type}\label{sl}
We will list up {\it five} types of isomonodromic deformations of $SL$-type. This part is the revision 
of the section 4.4 in \cite{Oka}. The isomonodromic deformation  of $SL$-type is 
\begin{equation*} \begin{aligned} 
\frac{\partial^2 u}{\partial x^2}&=p(x,t)u,\\
\frac{\partial u}{\partial t}=& A(x,t)\frac{\partial  u}{\partial x}
                             -\frac12\frac{\partial A(x,t)}{\partial x }u.
\end{aligned}\end{equation*}
The compatibility condition is given by
$$p_t(x,t)=2p(x,t)A_x(x,t)+A(x,t)p_x(x,t)-\frac12 A_{xxx}(x,t).$$
In the following list, $p(x,t)$ contains a Hamiltonian  $K$. The compatibility condition 
is  the Hamiltonian system with the Hamiltonian  $K$. Eliminating $z$, we obtain 
the Painlev\'e equation on $y$, which is an apparent singularity of 
the linear equation.  
We remark that R.~Fuchs studied the isomonodromic deformations of $SL$-type for P6 \cite{Fuchs1}.

\par\bigskip\noindent
{\it Type $(4),  (7/2)$}\ : P1\_2($\alpha, \beta$)\\
\begin{equation*}\begin{aligned} 
p(x,t)&= \alpha x^4 +4\beta x^3 + \alpha tx^2 +2\beta tx +2 K_{1,2}+\dfrac{3}{4(x-y)^2}-\dfrac{z}{x-y},\\
A(x,t)&= \dfrac{1}{2(x-y)},\\
K_{1,2} &=\dfrac{z^2}2-  \beta (2y^3+ ty) - \dfrac{\alpha}2 (y^4+ty^2).
\end{aligned} \end{equation*}

%% \begin{equation*}\begin{aligned} 
%% p(x,t)&= \alpha x^4 +4 x^3 + \alpha tx^2 +2tx +2 K_{1,2}+\dfrac{3}{4(x-y)^2}-\dfrac{z}{x-y},\\
%% A(x,t)&= \dfrac{1}{2(x-y)},\\
%% K_{1,2} &=\dfrac{z^2}2- 2y^3- ty- \dfrac{\alpha}2 (y^4+ty^2).
%% \end{aligned} \end{equation*}

\par\bigskip\noindent
{\it Type $(1)(3), (1)(5/2)$}\ : P4\_34${}^\prime(\alpha, \beta, \gamma$) for $\sigma=1$, 
                                 P4\_34($\alpha, \beta, \gamma$) for $\sigma=-1$.\\
\begin{equation*}\begin{aligned} 
p(x,t)&=\frac \gamma 2 \left( {x^2}+2 \sigma tx +  {  t^2 } \right)
+\frac \beta 2  ( x+ \sigma { t} ) +\frac{\alpha-1}{4x^2}+ 
\frac{\sigma K_{4,34}}{x}+ \frac{3}{4(x-y)^2} -\frac{ yz}{x(x-y)},\\
A(x,t)&=   \frac{\sigma x}{x-y},\\
K_{4,34} &= \sigma yz^2-\sigma z +\frac{\sigma(1-\alpha)}{4y}-\frac\beta2 y( \sigma y +{t })
                      -\frac{\sigma \gamma}2 y( y+\sigma t)^2.
\end{aligned} \end{equation*}

%% $$p(x,t)=\delta \left(\frac{x^2}{2} -tx +\frac {  t^2 }2\right)
%% +\ \frac x2 - \frac { t}2  +\frac{\alpha}{4x^2}+ 
%% \frac{K_{\rm IV}}{x}+ \frac{3}{4(x-y)^2}  -\frac{ z}{x(x-y)},$$
%% $$A(x,t)=  - \frac{  x}{x-y},$$
%% $$K_{\rm IV}  = yz^2-\frac{y^2}2 +\frac{ty}2-\frac{\alpha}{4y} 
%%                        -\delta \left(\frac{y^3}2 -ty^2+\frac{t^2y}{2} \right).$$

\par\bigskip\noindent
{\it Type $(2)^2, (2)(3/2), (3/2)^2$}\ :  P3${}^\prime$($\alpha,  \beta, \gamma, \delta$)\\
\begin{equation*}\begin{aligned} 
p(x,t)&=\frac{a_0t^2}{x^4}+\frac{a_0^\prime t }{x^3}+\frac{t K_{\rm III}^\prime}{ x^2}+
\frac{a_\infty^\prime}{ x } + a_\infty
+\frac{3}{4(x-y)^2}-\frac{ yz}{x(x-y)},\\
 A(x)&=   \frac{ yx}{t(x-y)},\\ 
 t K_{\rm III}^\prime &= y^2z^2 -yz -\frac{a_0 t^2}{ y^2}
-\frac{a_0^\prime t}{ y}- {a_\infty^\prime y} -a_\infty  {y^2},
\end{aligned} \end{equation*}
$$a_0= -\frac\delta{16},\ a_0^\prime=-\frac\beta8,\ 
a_\infty= \frac\gamma{16},\ a_\infty^\prime=\frac\alpha8.$$

\par\bigskip\noindent
{\it Type $(1)^2(2), (1)^2(3/2)$}\ :  P5($\alpha, \beta, \gamma, \delta$) 
\begin{equation*}\begin{aligned} 
p(x,t)&=\frac{a_1 t^2}{(x-1)^4}+\frac{K_{\rm V} t}{(x-1)^2 x}+\frac{a_2  t}{(x-1)^3}-\frac{z
   (y -1) y }{ x (x-1)(x-y )}+\frac{a_\infty}{(x-1)^2}+\frac{{a_0}}{x^2}+\frac{3}{4
   (x-y )^2},\\ 
 A(x)&=\frac{y-1}t \cdot \frac{x(x-1) }{ x-y },\\ 
 tK_{\rm V}&= y(y -1)^2   \left[-\frac{a_1 t^2}{(y-1)^4}-\frac{a_2 t}{(y -1)^3}+z^2-
   \left(\frac{1}{y }+\frac{1}{y -1}\right)z -\frac{a_\infty}{(y -1)^2}-\frac{a_0}{y 
   ^2}\right],
   \end{aligned} \end{equation*}
$$a_0= -\frac\beta 2-\frac14,\ a_1=-\frac\delta2,\ a_2=-\frac\gamma2, \ a_\infty=\frac12(\alpha+\beta)-\frac34.$$

\par\bigskip\noindent
{\it Type $  (1)^4$}\ :  P6($\alpha, \beta, \gamma, \delta$) 
\begin{equation*}\begin{aligned} 
p(x,t)&=\frac{a_0 }{x^2}+ \frac{a_1 }{(x-1)^2}+\frac{a_\infty}{x(x-1)} +\frac{b_1}{ (x-t)^2}
        +\frac{3}{4(x-y )^2} \\
        &\hskip 2cm + \frac{t(t-1) K_{\rm VI} }{x(x-1)(x -t)}-\frac{y(y-1)z}{x(x-1)(x -y)},\\ 
 A(x)&=\frac{y-t}{t(t-1)} \cdot \frac{x(x-1) }{ x-y },\\ 
K_{\rm VI}&= \frac{y(y -1)(y-t)}{ t(t-1)} \left[z^2 - \left(\frac{1}{y }+\frac{1}{y -1}\right)z
-\frac{a_0}{y^2} -\frac{a_1  }{(y-1)^2} -\frac{a_\infty}{y(y -1)} -\frac{b_1}{(y -t)^2} \right],
   \end{aligned} \end{equation*}
$$a_0=  -\frac\beta 2-\frac14,\ a_1= \frac\gamma2-\frac14,\  b_1= -\frac12\delta,\
a_\infty=\frac12(\alpha+\beta-\gamma+\delta-1).$$

\par\bigskip\noindent
If we set $\alpha=0$ in P1\_2, we obtain the standard isomonodromic deformations of $SL$-type for P1. 
If we set $\gamma=0$ in P4\_34, we obtain the  $SL$-type for P34. 
If we set $\gamma=0$ in P3${}^\prime$, we obtain  the $SL$-type for 
P3${}^\prime(D_7^{(1)})$.  If we set $\gamma=0, \delta=0$ in P3${}^\prime$, we obtain the  
$SL$-type for P3${}^\prime(D_8^{(1)})$.  
If we set $\delta=0$ in P5, we obtain the  $SL$-type for deg-P5.   

We show the standard isomonodromic deformations of $SL$-type for P2 and P4. 

\par\bigskip\noindent
P2($\alpha$):\\
\begin{equation*}\begin{aligned} 
p(x,t)&=x^4 +t x^2 +2\alpha x +2 K_{\rm II} +\frac{3}{4(x-y)^2}  -\frac{z}{x-y},\\
 A(x,t)&= \frac 12\cdot \frac{1}{x-y},\\
 K_{\rm II} & =\frac12 z^2 -\frac12 y^4 -\frac12 t y^2 -\alpha y. 
\end{aligned} \end{equation*}

\par\bigskip\noindent
P4($\alpha, \beta$):\\
\begin{equation*}\begin{aligned} 
p(x,t)&=\frac{{a_0}}{x^2}+\frac{K_{\rm IV}}{2x} +a_1 +\left(\frac{x+2t}4\right)^2 
+\frac{3}{4(x-y)^2}  -\frac{yz}{x(x-y)},\\
A(x)&=   \frac{2x}{x-y},\\ 
K_{\rm IV } & =2yz^2-2z -\frac{{2a_0}}{y} -2a_1 y-2y\left(\frac{y+2t}4\right)^2,
\end{aligned} \end{equation*}
$$a_0= -\frac\beta 8-\frac14,\ a_1=-\frac\alpha4.$$

\par\bigskip\noindent
{\it Remark.}\ In \cite{Oka}, there is a misprint in $K_{\rm IV}$.

%% P2  $y^{\prime\prime}= 2y^3+ty+\alpha $
%% $$y=\varepsilon \lambda + \varepsilon^{-5}, t =\varepsilon^2 t_1 -6 \varepsilon^{-6}$$
%% $$z= \varepsilon^{-1} \mu, K2= \frac12 z^2 -\frac12 y^4-\frac12t y^2- \alpha y$$

\end{document}